\documentclass[11pt]{amsart}
\usepackage{amssymb}
\allowdisplaybreaks[1]

\newtheorem{thm}{Theorem}[section]

\newtheorem{lemma}[thm]{Lemma}
\newtheorem{proposition}[thm]{Proposition}

\numberwithin{equation}{section}
\newcommand{\kP}{\mathbb{P}}
\newcommand{\kN}{\mathbb{N^{*}}}
\newcommand{\kkN}{\mathbb{N}}
\newcommand{\kZ}{\mathbb{Z}}
\newcommand{\kA}{\mathbb{A}}
\newcommand{\kB}{\mathbb{B}}
\newcommand{\kC}{\mathbb{C}}
\newcommand{\kR}{\mathbb{R}}
\newcommand{\kz}{\zeta}

\newcommand{\ksaut}{\vspace{0.4cm}\\ }

\newcommand{\kreel}[1]{\mathfrak{Re}(s)>{#1}}

\begin{document}
\title{Partial Euler products as a new approach to Riemann hypothesis}
\author{Jean-Paul JURZAK}
\maketitle

\begin{center}
{  {\small{Laboratoire Gevrey de Math\'ematique Physique, Universit\'e de Bourgogne, \\
Facult\'e des Sciences et Techniques\\
 BP 47870, 21078 Dijon Cedex\\
 e-mail : Jean-Paul.Jurzak@u-bourgogne.fr}}}
\end{center}

\vskip.5cm
\begin{abstract}
In this paper, we show that Riemann hypothesis (concerning zeros of the zeta function in the critical strip) is equivalent  to the analytic continuation of Euler products obtained by restricting the Euler zeta product to suitable subsets $M_k$, $k\geq 1$ of the set of prime numbers. Each of these Euler product  defines so a partial zeta function $ \zeta_{k}(s)$ equal to a Dirichlet series of the form $\sum\; \epsilon(n)/n^s$, with coefficients $ \epsilon(n)$ equal to  $0$ or $1$ as $n$ belongs or not to the population of integers generated by $M_k$. We show that  usual formulas of the arithmetic adapt themselves to such populations (Moebius, Mertens, Lambert series,...).  We envisage also the study of summations inside these populations and new functions (generalizations of the integer part function, of the harmonic series) directly connected to the existence of analytical continuations.
\end{abstract}
\vskip1.5cm
\section{Introduction}

The number $\pi(x)$ of prime numbers less or equal to $x$ is known to be equivalent to $\frac{x}{\ln(x)}$ and prime number theorem has been proved with the  
estimation:
$$ \pi(x) = \text{Li}(x) + O( x\; e^{-c\;\sqrt{\ln(x) }}) $$
for some suitable positive constant $c$:  an improvement (Vinogradov and Korobov) of the error term and  associated comments can be found in $[Nar]$ p.236. One importance
of Riemann hypothesis (which states that the non-real zeros of $\zeta(s)$ 
have real part equal to $\frac{1}{2}$) lies in the fact that the estimation 
$$ \pi(x) = \text{Li}(x) + O(\sqrt{x}\;\ln(x) ) $$
is probably true. Many arithmetical conjectures and various problems are directly connected to Riemann hypothesis. The zeta function has been studied intensively with developments in various directions. Literature on Riemann zeta function is known to be huge  and very diversified and this work hopes not to find  results already published.

One purpose of this work is to show that analytic continuation of partial zeta functions  
$$ \zeta_{k}(s)\;=\;\prod_{ p \in M_{k}} \frac{1}{( 1 - \frac{1}{p^s})} 
\quad \text{with} \quad M_{k} = \{p_{1},p_{k+1},p_{2k+1},p_{3k+1},\cdots\}$$
is equivalent to Riemann hypothesis: here, $k$ is an integer and
$\kP=(p_{n})_{n \in \kN}$ is the set of all primes. The obtaining of the analytical continuation of functions $\zeta_{k}(s)$ is in direct relation with the study of  subsets of $\kN$ generated by an arbitrary subset of $\kP$: this  approach is  a new facet of arithmetic. So, another purpose of this work is  to adapt classic formulas of  arithmetic to such subsets.

Every subset M, included in $\kP$ (the set of prime numbers) generates, by successive products  of the different elements of M, a subset of $\kN$, that is a sub-population of $\kN$.  As a consequence, we use  notation $pop(M)$ to indicate such a set, $pop$ term being the abbreviation of \textit{population}. When there is no ambiguity on the choice of $M$, one writes simply $pop$ in the place of $pop(M)$. By definition, $pop(\kP)=\kN$ and each function $ \zeta_{k}(s)$ is constructed by choosing for $M$ the  prime numbers taken by jumps of $k$ inside \textit{the set of index} of the  descriptive formula  $\kP=(p_{n})_{n \in \kN}$. Such an $M$ is called an arithmetical list of reason $k$ (see Conventions of section \ref{convent}): naturally, this does not mean that the $p_i$'s chosen constitute an arithmetical sub-sequence of $\kN$.

The integer part function $[x]$ of $x$, appropriate for the analysis in $\kN$, is then replaced  by the function
$$ N_{pop}(x) \; \equiv \; \text{card}([1;x]\cap \text{pop})$$
and another useful function is
$$  S_{pop}(x) \; \equiv \;\sum_{n \in \text{pop} \;,\; n \leq x}\; \frac{1}{n}$$ 
(generalization of the harmonic series)
 
The rigorous evaluation of $ N_{pop}(x)$ (or of $S_{pop}(x)$) insures the existence of the analytical continuation of partial zeta functions  
$ \zeta_{k}(s)$ to the open set 
$ \{ \mathfrak{Re}(s)>{\frac{1}{2}}\} - [ \frac{1}{2} ;1]$. More precisely, for $M$  an arithmetical list of reason $k$, these unproved evaluations are
$$ N_{\text{pop}(x)} = Cte \;\frac{x\;\ln(x) }{\sqrt[k]{\ln(x)}}  + \psi(x) $$
with $\psi(x)=O(x^{\frac{1}{2}+\epsilon})$ or
$$ S_{\text{pop}(x)} = Cte \;\sqrt[k]{\ln(x)} +  \phi(x) $$
with $\phi(x)= O( \frac{1}{x^{\frac{1}{2}-\epsilon}} )$. Formula \ref{fr:equilog} is  the  
natural equivalent of $ N_{pop}(x)$ and  we indicate in last section some consequences connected to this equivalence. In summary, when one passes through evaluation of $ N_{pop}(x)$ (to obtain the analytical continuation of $\zeta_{k}(s)$), Riemann's problem is returned to an arithmetical problem. One can also tempt others techniques to obtain  the existence of the analytical continuation of $\zeta_{k}(s)$.

\section{Preliminaries}

In the following, n stands for an integer, p stands for a prime number. We denote 
by  $\kP=(p_{n})_{n \in \kN}$ the set of all primes p (with $  p_1=2,p_2=3,p_3=5,p_4=7 \cdots $ ), by
$ \pi\;(x)$ the number of $p\in \kP $ satisfying $ 2 \leq p \leq x $. As usual, $\kz(s) = \sum_{n \in \kN} \frac{1}{n^s} $ is the Riemann zeta function, and
$Z_{\kz}$ the set of zeros of $\kz(s)$ in the half space $\kreel{\frac{1}{2}}$. It is a known fact that zeros of $\zeta(s)$ in the critical strip are displayed symetrically relatively to the line $\sigma = \frac{1}{2}$: thus, we may restrict some considerations to $\kreel{\frac{1}{2}}$.

We shall need in this paper, partial zeta functions $Z_{M}$ defined by an Euler product:
$$ Z_{M}(s)\;=\;\prod_{ p \in M} \frac{1}{( 1 - \frac{1}{p^s})} $$
where $M$ is a subset of $\kP$. Expanding each function:
$$ \frac{1}{( 1 - \frac{1}{p^s})}= 1 + \frac{1}{p^s} + \frac{1}{p^{2s}}+ \cdots $$
in power series, with $ p \in M $, and multiplying all expansions, we develop the Euler product as a Dirichlet series
$ \sum_{n \in pop_{M} } \frac{1}{n^s} $ , clearly analytic for $\kreel{1}$: here, n moves in the set $pop_{M}$ of all
integers of the form $ n = q_{1}^{\alpha_{1}} q_{2}^{\alpha_{2}}\cdots $ , with  all $\alpha_{i} \geq 0 \;$,$ \;\alpha_{i} \in \kN $ and all $q_{i}\in M$. The characteristic function of $pop_{M}$ will be denoted $ 1_{pop_{M}}(n)$, and equals 1 when $n \in pop_{M}$ and $0$ elsewhere. Thus, an equality of the form 
$$ \sum_{n \in pop_{M} } \frac{\alpha_{n}}{n^s} = \sum_{n \in \kN } 1_{pop_{M}}(n)\;\frac{\alpha_{n}}{n^s}
$$
means that the Dirichlet series considered on left side satisfies $\alpha_{n}=0$ whenever $ n \notin
 pop_{M}$. We also need:
$$  \eta_{M}(s)\;=\;\sum_{p \in M  } \frac{1}{p^s}$$
clearly analytic for $\kreel{1}$.
\begin{lemma} \label{le:ProD}
Let $ M \subset \kP$ and 
$$f(s)\;=\;\sum_{n \in pop_{M} } \frac{\alpha_{n}}{n^s}\; \qquad g(s)\;=\;\sum_{n \in pop_{M} } \frac{\beta_{n}}{n^s}$$ 
be two absolutely convergent Dirichlet series. Then, the Dirichlet product 
$$(fg)(s)\;=\;\sum_{n \in \kN } \frac{\gamma_{n}}{n^s} $$ 
of $f(s)$ and $g(s)$ satisfies $\gamma_{n}=0$ whenever $n \notin \text{pop}_M$.
\end{lemma}
The lemma follows from equality
$$\gamma_{n}= \sum_{kl=n ,\; k \geq1\;, l\geq1}\;\alpha_{k}\beta_{l}$$ 
and the fact that $k | n$ with $n\in pop_{M}$ implies $k \in pop_{M}$ and $l \in pop_{M}$.
\ksaut
\indent Let $Z_{M}(s)\;=\;\sum_{n \in pop_{M} } \frac{1}{n^s}$ be the associated partial zeta function (defined and analytic for $\kreel{1}$). One verifies that some of  known formulas
 concerning $\zeta(s)$ become now:
$$ \frac{Z_{M}(s-1) }{Z_{M}(s) } \;=\; \sum_{n \in pop_{M} }\;  \frac{\Phi{(n)}}{n^s} \qquad \text{for }\;\kreel{2}$$
$$  Z_{M}(s)\;Z_{M}(s-1)   \;=\; \sum_{n \in pop_{M} }\;  \frac{\sigma{(n)}}{n^s} \qquad \text{for }\;\kreel{2}$$
$$  Z^{2}_{M}(s)   \;=\; \sum_{n \in pop_{M} }\;  \frac{\tau{(n)}}{n^s} \qquad \text{for }\;\kreel{1}$$
$$ \frac{1 }{Z_{M}(s) } \;=\; \sum_{n \in pop_{M} }\;  \frac{\mu{(n)}}{n^s} \qquad \text{for }\;\kreel{1}$$
$$  (1+\eta_{M}(s))\;Z_{M}(s)   \;=\; \sum_{n \in pop_{M} }\;  \frac{\nu{(n)}}{n^s} \qquad \text{for }\;\kreel{1}$$
$$
\sum_{n \in pop_{M} }\;  \frac{\ln{(n)}}{n^s}
\;=\; Z_{M}(s)   \;\sum_{p \in M }\;  \frac{\ln{(p)}}{p^s-1} \qquad \text{for }\;\kreel{1}$$
where, as usual:\\
$\Phi(n)$ = the Euler Phi-function; \\
$\sigma(n) = \Sigma_{i | n}\;i$, the sum of the divisors of n;\\
$\tau(n)= \Sigma_{i | n}\;1$ , the number of  divisors of n;\\
$\nu(n)=\Sigma_{p | n}\;1$, the number of distinct prime factors of n;\\
$\mu(n)$ the Moebius function,$\;\mu(1)=1$ , $\mu(n)=0$ if n is divisible by a square $> 1$, $ = (-1)^{ \nu(n)}$ in others cases.
\ksaut
\indent Subsets $M$ of main interest will be {\it{arithmetical lists}} obtained from the list of all primes:
$$ \kP = \{2,3,5, 7, 11, 13, 17, 19, 23, 29, 31, 37, 41, 43, 47, 53, 59, 61, 67, 71, 
  73, 79, 83 ,\cdots \}$$
\textbf{Conventions:} \label{convent} \\
\textsl{An arithmetical list of $\kP$ having a reason $ r >0$ is a subset $M \subset \kP$ of the form 
$$M=\{p_{r_{0}},p_{r_{0}+r},p_{r_{0}+2r},\cdots\}$$
for some $r_{0}>0$ with $r_{0}\leq r$.\\
Whenever a function $f_{M}(\cdots)$ is defined or depends on an
arithmetical list
 determined by its two first elements $a \in M$ and $b \in M$ with $a=p_{r_{0}}$ and 
 $b=p_{r_{0}+r}$, we shall write
$$ M\;=\; M_{a,b}$$
$$ f_{M}(\cdots)\;=\; f_{a,b}(\cdots)$$
Moreover, when the arithmetical list $M_{a,b}$ starts with the first prime (ie one 
has  $a=p_{1}=2$), we write simply
$$M_{a,b}\;=\;M_{r}$$
$$f_{a,b}(\cdots)\;=\;f_{r}(\cdots) $$
where $r$ is the reason of the list.}
\ksaut
\indent For example, 
$ M = \{5, 13,  23, 37, 47,  61, 73,\cdots \}$
is an arithmetical list with first element $p_{3}=5$ and reason $r=3$.

We observe that many results concerning functions built with an arithmetical list are easily adapted to  functions built with subsets $M \subset \kP$ which are finite union
of arithmetical lists having all the same reason $ r >0$ (such an $M$ has the
following property:  $p_{k} \in M$ implies $p_{k+r} \in M$ and
$p_{k-r} \in M$ whenever $k-r\geq1$).\\
With our conventions, choosing an $M$ of the form $M_{a,b}$, one may write 
$Z_{M}(s)= Z_{a,b}(s)$ and $ 1_{pop_{M}}(n)= 1_{pop_{a,b}}(n)$. In fact,
$$ \zeta_{k}(s)\;=\;\prod_{ p \in M_{k}} \frac{1}{( 1 - \frac{1}{p^s})} 
\quad \text{with} \quad M_{k} = \{p_{1},p_{k+1},p_{2k+1},p_{3k+1},\cdots\}$$
$$ \zeta_{2}(s) \;=\;\prod_{ n \in \kkN \; n\geq 0} \frac{1}{( 1 - \frac{1}{p_{2n+1}^s})} \;=\; 
\frac{1}{( 1 - \frac{1}{2^s})}\;
\frac{1}{( 1 - \frac{1}{5^s})}\;
\frac{1}{( 1 - \frac{1}{11^s})}\;
\cdots
\frac{1}{( 1 - \frac{1}{p_{2n+1}^s})}
\cdots
\;=\;Z_{2,5}(s) 
$$
$$ \zeta_{3}(s) \;=\;\prod_{ n \in \kkN \; n\geq 0} \frac{1}{( 1 - \frac{1}{p_{3n+1}^s})} \;=\; 
\frac{1}{( 1 - \frac{1}{2^s})}\;
\frac{1}{( 1 - \frac{1}{7^s})}\;
\frac{1}{( 1 - \frac{1}{17^s})}\;
\cdots
\frac{1}{( 1 - \frac{1}{p_{3n+1}^s})}
\cdots
\;=\;Z_{2,7}(s) 
$$
$$ Z_{3,7}(s) \;=\;\prod_{ n \in \kkN \; n\geq 0} \frac{1}{( 1 - \frac{1}{p_{2n+2}^s})} \;=\; 
\frac{1}{( 1 - \frac{1}{3^s})}\;
\frac{1}{( 1 - \frac{1}{7^s})}\;
\frac{1}{( 1 - \frac{1}{13^s})}\;
\cdots
\frac{1}{( 1 - \frac{1}{p_{2n+2}^s})}
\cdots 
$$
$$ Z_{2,5}(s)  
=
{1} +\frac{1}{ 2^s} +\frac{1}{ 4^s} +\frac{1}{ 5^s} +\frac{1}{ 8^s} +\frac{1}{ 10^s}
+\frac{1}{ 11^s} +\frac{1}{ 16^s} +\frac{1}{ 17^s} +\frac{1}{ 20^s}  +\frac{1}{ 22^s}
+\frac{1}{ 23^s} +\frac{1}{ 25^s} +\frac{1}{ 31^s} +\frac{1}{ 32^s}  + \cdots $$
$$ Z_{3,7}(s) = 
1 + \frac{1}{3^{s}} + \frac{1}{7^{s}} + \frac{1}{9^{s}} + \frac{1}{13^{s}} + \frac{1}{19^{s}} + 
\frac{1}{21^{s}} + \frac{1}{27^{s}} + \frac{1}{29^{s}} + \frac{1}{37^{s}} + \frac{1}{39^{s}} + 
  \frac{1}{43^{s}} + \frac{1}{49^{s}}+\cdots
$$
$$\eta_{3,7}(s) \;=\; \frac{1}{ 3^s }+\frac{1}{ 7^s }+\frac{1}{ 13^s }+\frac{1}{ 19^s }+\frac{1}{ 29^s }+\frac{1}{ 37^s }
+\frac{1}{ 43^s }+\cdots
$$
$$\eta_{2}(s) \;=\; \frac{1}{ 2^s }+\frac{1}{ 5^s }+\frac{1}{ 11^s
}+\frac{1}{ 17^s }+\frac{1}{ 
  23^s }+\frac{1}{ 31^s }+\frac{1}{ 41^s}+\frac{1}{ 47^s}+\cdots
\;=\; \eta_{2,5}(s)
$$
$$\eta_{1}(s) = \sum_{p \in \kP} \frac{1}{p^s} \equiv \frac{1}{ 2^s }+\frac{1}{ 3^s }+\frac{1}{ 5^s }+\frac{1}{ 7^s }+\frac{1}{ 11^s
}+\frac{1}{ 13^s }+\frac{1}{ 17^s }+\frac{1}{ 19^s }+\frac{1}{ 
  23^s }+\frac{1}{ 29^s }+\frac{1}{ 31^s }+\cdots=\eta_{2,3}(s)
$$
\ksaut
\indent A central result, proved in section 5, is:
\begin{thm} \label{th:RieEqu}
Let $\Omega_{\frac{1}{2}}$ be the open set $\Omega_{\frac{1}{2}} \;=\;\{ \kreel{\frac{1}{2}} \;;\; s \notin [\frac{1}{2};1]\; \}$. 
The following properties are equivalent:
\begin{enumerate}
\item The Riemann hypothesis holds.
\item For every $i \geq 2 $ , the function $\zeta_{i}(s)$ extends as an analytic (resp. meromorphic) function in $\Omega_{\frac{1}{2}}$.
\end{enumerate} 
\end{thm}

\section{Basic Properties}

\begin{lemma} \label{le:ProEqu}
Let ($\alpha_{n}$) and ($\beta_{n}$) be two increasing sequences 
of reals such that $\alpha_{n}\leq \beta_{n} \leq \alpha_{n+1} $ for all $n$, 
$ \alpha_{1} > 1$ and 
$ \Sigma_{n}\;\frac{1}{\alpha_{n}^{\lambda}} < +\infty$ for every $\lambda > 1$.\\ 
Then, the following properties holds:
\begin{enumerate}
\item Euler products
$$  
 E_{1}(s) = \prod_{ n } \frac{1}{( 1 - \frac{1}{\alpha_{n}^s})} 
 \qquad
E_{2}(s) = \prod_{ n } \frac{1}{( 1 - \frac{1}{\beta_{n}^s})} 
$$
are well defined and analytic for $\kreel{1} $.
\item There exists a function $f(s)$, defined and analytic for $\kreel{0} $,
which never vanishes on this set,
such that $E_{1}(s) \;=\; E_{2}(s)  f(s)$ when $\kreel{1} $.
\end{enumerate}
\end{lemma}
If one knows besides that $  E_{1}(s) $ has a meromorphic (resp. analytic) 
continuation for $\kreel{a} $ for some $0<a<1$ (for example, $E_1$ may be the Riemann zeta function as in $[Gr$-$Sc]$),  our lemma implies that $E_2$ has also a meromorphic (resp. analytic) continuation for $\kreel{a} $. Therefore, functions $  E_{1} $ and $  E_{2} $
have the same zeros and same poles for $\kreel{a} $. 
\ksaut
{\indent \it{Proof:}}
It is known that $ \big|(1 - (1-z)e^z)\big| \;\leq \; \big| z \big|^{2}$ for $ \big|z\big|\leq 1$, 
hence the infinite product $\prod_{ n }\;( 1 - \frac{1}{\alpha_{n}^s})\; e^{\frac{1}{\alpha_{n}^s}}$ defines of an analytic function for $\kreel{\frac{1}{2}} $. The equation
$( 1 - \frac{1}{\alpha_{n}^s})\;=\;0$
has no zeros in variable $s$, hence the infinite product never vanishes. 
Clearly,
$$ E_{1}(s) = \; e^{\Sigma_{n}\;\frac{1}{\alpha_{n}^s}}\;\prod_{ n } \frac{1}{( 1 - \frac{1}{\alpha_{n}^s})\; e^{\frac{1}{\alpha_{n}^s}} } $$
for $\kreel{1} $ and this formula clearly defines an analytic function on that set. Since $ \Sigma_{n}\;\frac{1}{\beta_{n}^{\lambda}} < +\infty $ for every $\lambda > 1$, we get same conclusions for
$\prod_{ n }\;( 1 - \frac{1}{\beta_{n}^s})\; e^{\frac{1}{\beta_{n}^s}}$: this proves 1).\\ 
We now show that there exists a function $f(s)$, defined and analytic for $\kreel{\frac{1}{2}} $,
which satisfy  formula of 2). Assuming $ n \geq 1$, one has, for $\kreel{1}$:
$$ E_{1}(s) = E_{2}(s) \;  \frac
{\prod_{ n }( 1 - \frac{1}{\beta_{n}^s})\; e^{\frac{1}{\beta_{n}^s}} }
{\prod_{ n }( 1 - \frac{1}{\alpha_{n}^s})\; e^{\frac{1}{\alpha_{n}^s}}}
\;e^{ \frac{1}{\alpha_{1}^{s}} -\frac{1}{\beta_{1}^{s}} + \cdots +
 \frac{1}{\alpha_{n}^{s}} -\frac{1}{\beta_{n}^{s}}+ \cdots
 }$$
We put:
$$ f(s)\;=\;\frac
{\prod_{ n }( 1 - \frac{1}{\beta_{n}^s})\; e^{\frac{1}{\beta_{n}^s}} }
{\prod_{ n }( 1 - \frac{1}{\alpha_{n}^s})\; e^{\frac{1}{\alpha_{n}^s}}}
\;e^{ \frac{1}{\alpha_{1}^{s}} -\frac{1}{\beta_{1}^{s}} + \cdots +
 \frac{1}{\alpha_{n}^{s}} -\frac{1}{\beta_{n}^{s}}+ \cdots
 }$$
From Dirichlet series,
 we see that the series 
 $\frac{1}{\alpha_{1}^{s}} -\frac{1}{\beta_{1}^{s}} + \cdots +
 \frac{1}{\alpha_{n}^{s}} -\frac{1}{\beta_{n}^{s}}+ \cdots
$
is analytic in the half-plane $\kreel{0}$. Thus, $f(s)$ has the required properties for
$\kreel{\frac{1}{2}} $.\\
In the general case, taking $k \in \kN$, we use the Weierstrass factor 
$$W_{k}(z)=(1-z)\;e^{z+ \frac{z^2}{2}+\cdots+\frac{z^k}{k}}$$
 having the property
$ \big|1 -W_{k}(z) \big| \;\leq \; \big| z \big|^{k+1}$ for $ \big|z\big|\leq 1$.
Writting:
$$ E_{1}(s) = E_{2}(s) \;  \frac
{\prod_{ n } W_{k}(\frac{1}{\beta_{n}^s})}
{\prod_{ n } W_{k}(\frac{1}{\alpha_{n}^s})}
\;e^{ \frac{1}{\alpha_{1}^{s}} -\frac{1}{\beta_{1}^{s}} + \cdots +
 \frac{1}{\alpha_{n}^{s}} -\frac{1}{\beta_{n}^{s}}+ \cdots
 }
\cdots
\;e^{ \frac{1}{k}\; \big(
\frac{1}{\alpha_{1}^{ks}} -\frac{1}{\beta_{1}^{ks}} + \cdots +
 \frac{1}{\alpha_{n}^{ks}} -\frac{1}{\beta_{n}^{ks}}+ \cdots
\big)
 }
$$
our conclusion follows for $\kreel{\frac{1}{k+1}} $. Now, $f(s)$ is the analytic continuation of $\frac{E_{1}(s)}{E_{2}(s)}$ for $\kreel{\frac{1}{k+1}} $ and all $k \in \kN$: this proves the lemma.
\begin{proposition} \label{pr:ZetZek}
Let $k \geq 2,\; k\in \kN$ and $\Omega$ be a simply connected open set contained 
in $ \kreel{0} $ in which the Riemann function $\zeta(s)$ never vanishes. The following properties holds:
\begin{enumerate}
\item  $\zeta_{k}(s)$ extends analytically for $s \in \Omega $.
\item There exists a function $g_{k}(s)$ defined and analytic for $\kreel{0} $,
having no zeros on this half-plane,
such that 
\begin{equation} \label{fr:zetgzetk}
\zeta(s) \;=\; g_{k}(s)\;\zeta_{k}(s)^{k}   \quad \text{for} \quad s\in \Omega
\end{equation}
\end{enumerate}
\end{proposition}
It follows that 
\begin{equation} \label{form:akzetak}
\zeta_{k}(s) \;\sim\; \frac{A_{k}}{\sqrt[k]{1-s}}\quad \text{as}\quad s \rightarrow 1
\end{equation}
for some constant $A_{k}$.
\ksaut
{\indent \it{Proof:}}
One has $ \zeta_{k}(s)\;=\;\prod_{ p \in M_{k}} \frac{1}{( 1 - \frac{1}{p^s})} $.
 For $i \geq1 , i\leq k$, we define 
$$A_{i}= \{p_{i},p_{k+i},p_{2k+i},p_{3k+i},\cdots\}
\quad \text{and} \quad
 Z_{p_{i},p_{k+i} }(s)\;=\;\prod_{ p \in A_{i}} \frac{1}{( 1 - \frac{1}{p^s})}$$
Each $A_{i}$ is an arithmetical  list of $\kP$ with $A_{1}=M_{k}$ and $\kP = \bigcup_{i=1}^{i=k}A_{i}  $.
One has 
$ \zeta(s)\;=\; \prod_{ i=1}^{i=k}\;Z_{p_{i},p_{k+i} }(s) $ and $\zeta_{k}(s)=Z_{p_{1},p_{k+1} }(s)$.
 Using lemma \ref{le:ProEqu}, we get for every $i>1$ an analytic function $f_{i}(s)$   having no zeros 
for $\kreel{0} $ such that  $Z_{p_{i},p_{k+i} }(s)
=f_{i}(s)\zeta_{k}(s)$ for $\kreel{1}$.
Putting $g_{k}(s)=f_{1}(s)f_{2}(s) \cdots f_{k-1}(s)$, we get
$\zeta(s) \;=\; g_{k}(s)\;\zeta_{k}(s)^{k} $ for  $\kreel{1} $.\\
Let us now assume $s \in \Omega$. There exists an analytic function $h(s)$ analytic in $\Omega$ such that 
$\frac{\zeta(s)}{ g_{k}(s)}=e^{h(s)}$ for $s \in \Omega$. Thus 
$$ \bigg(\frac{ \zeta_{k}(s)} {e^{\frac{h(s)}{k}}}\bigg)^{k}\;=\;1 \quad \text{for} \quad \kreel{1} $$
therefore, for some fixed $\alpha \in \kZ$,
$ \zeta_{k}(s) \;=\; e^ {2 i \frac{\alpha}{k} \pi} \;e^{ \frac{ h(s)}{k} }
\quad \text{for} \quad \kreel{1}$. This formula defines now the analytic continuation of $ \zeta_{k}(s)$ for $s \in \Omega$
and formula \ref{fr:zetgzetk}
holds for $s \in \Omega$ by uniqueness of analytic continuations.
\ksaut
\indent The following lemma is contained  in the proof of lemma \ref{le:ProEqu}.
\begin{lemma} \label{le:ZetWexp}
Let $M \subset \kP$. The function 
$ W_{M}(s) = \prod_{  p \in M } 
\frac{1}{
( 1 - \frac{1}{p^s}) e^{\frac{1}{p^s}}} $
is analytic for $ \kreel{\frac{1}{2}}$ and has no zeros on this open set. It satisfies:
\begin{equation}
 Z_{M}(s)  = W_{M}(s) \; e^{\eta_{M}(s)} \quad \text{for} \quad  \kreel{1}
\end{equation}
\end{lemma}
$\;$

\section{Functions $\eta_{M}(s)$}

As defined previously, 
$ \displaystyle \eta_{M}(s)\;=\;\sum_{p \in M} \frac{1}{p^s}$
and
$\eta_{1}(s) = \sum_{p \in \kP} \frac{1}{p^s} =\eta_{2,3}(s)
$.
\begin{proposition} \label{pr:EtaEtM}
Let $M=\{p_{r_{0}},p_{r_{0}+r},p_{r_{0}+2r},\cdots\} $ be an arithmetical list of $\kP$ with a reason $r>1$ and $\Omega$ be an open  connected  set contained in $\kreel{0}$. Then, one has:
\begin{equation}  \label{fr:retam}
 \eta_{1}(s) = r\;\eta_{M}(s) + w(s)\quad \text{for}\quad \kreel{1}
\end{equation}
where $w(s)$ is analytic for $\kreel{0}$.\\
Therefore, 
$\eta_{M}(s)$ has an analytic (resp. meromorphic) continuation to $\Omega$ if and only if $\eta_{1}(s)$ has the same property.
\end{proposition}
{\indent \it{Proof:}}
Let $M_{r}=\{p_{1},p_{1+r},p_{1+2r},\cdots,p_{1+nr},\cdots\}\subset \kP$. 
 For every $0 \leq j<r$, we put:
$$M^{(j)}=\{p_{1+j},p_{1+j+r},p_{1+j+2r},\cdots,p_{1+j+nr},\cdots\}$$
to be the shifted list of $M_{r}$, hence $M^{(0)}=M_{r}$. One has
$$\eta_{M^{(0)}}(s)-\eta_{M^{(j)}}(s) \equiv \sum_{p \in M_{r}} \frac{1}{p^s}-\sum_{p \in M^{(j)}} \frac{1}{p^s}$$
for $\kreel{1}$, thus
$$\eta_{M^{(0)}}(s)-\eta_{M^{(j)}}(s)
= \frac{1}{p_{1}^s}-\frac{1}{p_{1+j}^s} + \frac{1}{p_{1+r}^s}-\frac{1}{p_{1+j+r}^s}
+\cdots+ \frac{1}{p_{1+nr}^s}-\frac{1}{p_{1+j+nr}^s}+\cdots
$$
The Dirichlet series on the right hand side is an analytic function $w_{j}(s)$ in $\kreel{0}$ for every $0 \leq j<r$. From
$$ \eta_{1}(s)= \eta_{M^{(0)}}(s)+\eta_{M^{(1)}}(s)+\cdots+\eta_{M^{(r-1)}}(s)$$
we get
$$ \eta_{1}(s)= \eta_{M^{(0)}}(s)+ \big( \eta_{M^{(0)}}(s)-w_{0}(s) \big)
+\cdots+\big( \eta_{M^{(0)}}(s)-w_{r-1}(s) \big)$$
which leads to
$$ \eta_{1}(s)= r \;\eta_{M^{(0)}}(s)-W(s)$$
with $W(s) \equiv w_{0}(s)+\cdots+w_{r-1}(s)$. This proves \ref{fr:retam} for $M=M^{(0)}$.\\
The arithmetical sublist 
$M=\{p_{r_{0}},p_{r_{0}+r},p_{r_{0}+2r},\cdots\} $  is one of the $M^{(j)}$ (with $j= r_{0}-1$) and relation
$\eta_{M^{(0)}}(s)-\eta_{M^{(j)}}(s) = w_{j}(s)$ implies formula \ref{fr:retam}. It follows
 from this formula that  $\eta_{M}(s)$ extends as an analytic
(resp. meromorphic) function  to $\Omega$ if and only if $\eta_{1}(s)$ has the same property. 
\ksaut
\indent We recall that $Z_{\kz}$ stands for the set of zeros of $\kz(s)$ in the half space $\kreel{\frac{1}{2}}$.
\begin{proposition} \label{pr:DEta}
Let $M$ be an arithmetical list.
\begin{enumerate}
\item The derivative
$  \quad \frac{d}{ds} \eta_M(s)\;=\;-\sum_{p \in M} \; \frac{\ln(p)}{p^s} \quad $
extends as a meromorphic function in $\kreel{0}$
with simple poles for singular points.
\item $\eta_{M}(s)$ can be continued to an analytic function
in the open simply connected 
$$ \Omega_{c}\;=\; \{\kreel{0}\} - \bigcup_{\alpha \in Z_{\zeta}}[0;\alpha] - [0;1] $$ 
\end{enumerate}
\end{proposition}
Assertion 2) is more precise than [Tit] page 182.
\ksaut
{\indent \it{Proof:}}
Let $r>0$ be the reason of the list. By proposition \ref{pr:EtaEtM},
$ \eta_{1}(s) = r\;\eta_{M}(s) + w(s)$ for $\kreel{1}$
, where $w(s)$ is analytic for $\kreel{0}$. Taking the derivative, we are reduced
to the case $\eta_{M}=\eta_{1}$ for which this property is true  for $\kreel{0}$. This 
proves 1) since our formula now defines the holomorphic continuation.\\
For 2), we note that $\Omega_{c}$ has the property that $s \in \Omega_{c}$
and $\lambda \geq1$ implies that $\lambda s \in\Omega_{c}$. 
Since the Riemann function $\zeta(s)$ never vanishes on $\Omega_{c}$, it follows from proposition \ref{pr:ZetZek} that
$$\zeta(s) \;=\; g_{k}(s)\;Z_{M}(s)^{k}   \quad \text{for} \quad s\in \Omega_{c}$$
for some suitable analytic function $g_{k}$ defined on $\Omega_{c}$ which never vanishes on this set. Now, the function:
$$ W_{M}(s) = \prod_{  p \in M } 
\frac{1}{
( 1 - \frac{1}{p^s}) e^{\frac{1}{p^s}}} $$
is analytic for $ \kreel{\frac{1}{2}}$ , has no zeros on this half-plane and satisfies $ Z_{M}(s)  = W_{M}(s) \; e^{\eta_{M}(s)}$ for $\kreel{1}$. Thus,
$ \frac{ Z_{M}(s)}{ W_{M}(s)} \;=\; e^{\eta_{M}(s)}$ 
has no zeros in $\Omega_{c}$ and
this leads to the analytic continuation on $\Omega_{c} \cap \{ \kreel{\frac{1}{p+1}} \}$
with $p=1$. We now iterate this process as follows. For $p\in \kP$ with $p \geq1$
$$ W_{M}^{(p)}(s) = \prod_{  p \in M } 
\frac{1}{
( 1 - \frac{1}{p^s}) e^{\frac{1}{p^s} + \frac{1}{2}\frac{1}{p^{2s}} +\cdots+
\frac{1}{p}\frac{1}{p^{2s}} }} $$
is analytic for $ \kreel{\frac{1}{p+1}}$ , has no zeros on this half-plane and satisfies:
$$ Z_{M}(s)  = W_{M}^{(p)}(s) \; e^{\eta_{M}(s)+ \frac{1}{2}\eta_{M}(2s) +\cdots+ \frac{1}{p}\eta_{M}(ps)} \quad \text{for} \quad  \kreel{1}$$
It follows that $\eta_{M}(s)+ \frac{1}{2}\eta_{M}(2s) +\cdots+ \frac{1}{p}\eta_{M}(ps)$
has an analytic continuation $f_{p+1}(s)$ on $\Omega_{c} \cap \{ \kreel{\frac{1}{p+1}} \}$. By induction, the analytic continuation $\overset{p}{\eta}_{M}(s)$ of $\eta_{M}(s)$ being obtained on $\Omega_{c} \cap \{ \kreel{\frac{1}{p}} \}$, one has, for $\kreel{\frac{1}{p+1}}$
$$ f_{p+1}(s) \;=\; \eta_{M}(s)+ \frac{1}{2}\overset{p}{\eta}_{M}(2s)+\cdots
+ \frac{1}{p}\overset{p}{\eta}_{M}(ps)$$
since $ \mathfrak{Re}(2s) > \frac{1}{p},\cdots,\mathfrak{Re}(ps) > \frac{1}{p}$. Thus,
$$ f_{p+1}(s) - \frac{1}{2}\overset{p}{\eta}_{M}(2s)-\cdots
- \frac{1}{2}\overset{p}{\eta}_{M}(ps)$$
defines the analytic continuation $\overset{p+1}{\eta}_{M}(s)$ of $\eta_{M}(s)$
to the open half-plane $\kreel{\frac{1}{p+1}}$.
\ksaut
\indent The following proposition is known for $M = \kP$.
\begin{proposition} \label{pr:EtaLog}
Let $M \subset \kP$.
\begin{enumerate}
\item  one has, for $\kreel{1}$, 
$$ \sum_{p \in M} \; \frac{1}{p^s}\;=\; \sum_{n=1}^{+\infty} \; \frac{\mu(n)}{n}\; \ln(Z_{M}(ns) )$$
\item Let $l_{\text{pop}_{M}}(x)\;=\; \sum_{k \in \text{pop}_{M}}\;\frac{x^k}{k}$ for $|x| <1$. One has:
$$ \sum_{n \in \text{pop}_M} 
\frac{\mu(n)}{n}\;l_{\text{pop}_{M}}(x^n)\;
=\;x$$
$$ l_{\text{pop}_{M}}(x)\;=\; \sum_{n \in \text{pop}_{(\kP-M)}} 
\frac{\mu(n)}{n}\;\ln(\frac{1}{1-x^n})
$$
It follows that
$$ \sum_{n \in \text{pop}_M} 
\frac{\mu(n)}{n}\;\bigg(\sum_{p \in M}\;l_{\text{pop}_{M}}(\frac{1}{p^{ns}})\bigg)\;
=\;\sum_{p \in M}\frac{1}{p^{s}}$$
\end{enumerate}
\end{proposition}
The function $l_{\text{pop}_{M}}(x)$ coincide with $\ln(\frac{1}{1-x})$ when $M=\kP$ (because $\text{pop}_{M}=\kN$ in this case). 
Taking  $\text{pop}_{M} =\text{pop}_{2,5}$ in 2), one has:
$$l_{pop}(x)= x + \;\frac{x^2}{2}+ \;\frac{x^4}{4}
+ \;\frac{x^5}{5}
+ \;\frac{x^8}{8}
+\cdots$$
and our first summation is
$$ 
\;\frac{\mu(1)}{1}\;l_{pop}(x)\;+
\;\frac{\mu(2)}{2}\;l_{pop}(x^2)\;+
\;\frac{\mu(4)}{4}\;l_{pop}(x^4)\;+
\;\frac{\mu(5)}{5}\;l_{pop}(x^5)\;+
\;\frac{\mu(8)}{8}\;l_{pop}(x^8)\;+
\cdots
$$
\ksaut
{\indent \it{Proof:}}
For 1), we use formula
$ \sum_{n=1}^{+\infty} \; \frac{\mu(n)}{n}\; \ln(\frac{1}{1-x^n})\;=\;x $
valid for $|x|<1$. Choosing $x = \frac{1}{p^s}$ for $p\in M$ and $\kreel{1}$, it remains to sum expansions for all $p $.\\
For 2), one has
$$ \sum_{n \in \text{pop}_M} 
\frac{\mu(n)}{n}\;l_{\text{pop}_{M}}(x^n)\;
=\; \sum_{n \in \text{pop}_M} 
\frac{\mu(n)}{n}\;\sum_{k \in \text{pop}_{M}}\;\frac{x^{nk}}{nk}
\;=\; \sum_{n \in \text{pop}_M} \frac{x^{n}}{n}
\;\sum_{\substack{ k\;\mid\; n \\ k \in \text{pop}_{M}}}\;\frac{\mu(\frac{n}{k})}{\frac{n}{k}}$$
Since conditions $k\mid n $ with $ k \in \text{pop}_{M}$
 and $k\mid n $ are identical when $n \in \text{pop}_{M}$
, using for $n >1$
 $$\sum_{d|n}\frac{\mu(d)}{d \frac{n}{d}}=\frac{1}{n}\;\sum_{d|n}\mu(d)=0$$
one has, from $\mu(1)=1$ 
$$ \sum_{n \in \text{pop}} 
\frac{\mu(n)}{n}\;l_{pop}(x^n)\;
=\;x$$
This proves first formula of 2). Formula
$$ l_{\text{pop}_{M}}(x)\;=\; \sum_{n \in \text{pop}_{(\kP-M)}} 
\frac{\mu(n)}{n}\;\ln(\frac{1}{1-x^n})
$$
follows from inclusion-exclusion principle. Last formula is obtained
with  $x= \frac{1}{p^{s}}$ in first formula of 2) and summation for all $p\in M$.

\section{Equivalences between analytic continuations and Riemann hypothesis}

As seen previously, $\Omega_{\frac{1}{2}} = \{ \kreel{\frac{1}{2}}\} - [ \frac{1}{2} ;1]$.
\begin{proposition} \label{pr:Rie5Equ}
The following properties are equivalent: 
\begin{enumerate}
\item $\eta_{1}(s)$ extends as an analytic function in $\Omega_{\frac{1}{2}}$. 
\item $\eta_{1}(s)$ extends as a meromorphic function in $\Omega_{\frac{1}{2}}$.
\item The Riemann function $\zeta(s)$ never vanishes in $\Omega_{\frac{1}{2}}$.
\item For one $i \in \kN\;,\; i\geq2$
(resp. for all $ i \in \kN$), the function $\zeta_{i}(s)$ extends as an analytic function  in   $\Omega_{\frac{1}{2}}$ and never vanishes in that open set.
\item $\eta_{1}^{n}(s)$ extends as a meromorphic function in $\Omega_{\frac{1}{2}}$, for one $n \in \kN\;,\; n\geq1$ (resp. for all $ n \in \kN$). 
\end{enumerate}
\end{proposition}
{\indent \it{Proof:}}
Let $M$ be an arithmetical list of $\kP$. From proposition \ref{pr:EtaEtM}, 
$\eta_{M}(s)$ has an analytic (resp. meromorphic) continuation to $\Omega$ if and only if $\eta_{1}(s)$ has the same property: and one has
$ Z_{M}(s)  = W_{M}(s) \; e^{\eta_{M}(s)}$ for $\kreel{1}$.
where $W_{M}$ is an analytic function for $ \kreel{\frac{1}{2}}$ and has no zeros on this half plane.\\
Assume that $1)$ is satisfied. Taking $M=M_{i}$ with $i\geq1$, one has $ \zeta_{i}(s)  = W_{M_{i}}(s) \; e^{\eta_{i}(s)}$ for $\kreel{1}$, 
thus the right side hand has an analytic continuation to $\Omega_{\frac{1}{2}}$: thus,
 $1\Rightarrow 3$ and $1\Rightarrow 4$. Conversely, if 4) is satisfied for one $i\geq2$, then
$ \displaystyle \frac{\zeta_{i}(s) }{W_{M_{i}}(s) } = e^{h(s)}$
for some analytic function $h=h(s)$ defined on $\Omega_{\frac{1}{2}}$. It follows that
$\eta_{i}(s)-h(s)$ is constant for $\kreel{1}$, hence $\eta_{i}(s)$ has an analytic continuation in $\Omega_{\frac{1}{2}}$. Hence,  $3\Rightarrow 1$ and 
$4\Rightarrow 1$.\\
Let us assume 2). Let $\alpha$ be a pole of $\eta(s)$ in $\Omega_{\frac{1}{2}}$. Formula
$ \zeta(s)  = W(s) \; e^{\eta(s)}$ has an analytic continuation in $\Omega_{\frac{1}{2}}-
\{\text{the set of poles of}\;\eta\}$. Taking a small neighborhood $V(\alpha)$ of $\alpha$ not containing $\alpha$, we see that the right (resp. left) hand side our the formula is unbounded (resp. bounded) on $V(\alpha)$ leading to a contradiction: hence $1\Leftrightarrow 2$.\\
 Let us assume 5) for some $ n\geq2$. The derivative 
$\displaystyle \frac{d}{ds}\eta^n(s)=n \eta^{n-1}(s)
\frac{d}{ds}\eta(s)$ is meromorphic in $\Omega_{\frac{1}{2}}$. Since $\frac{d}{ds}\eta(s)$ is meromorphic in $\kreel{0}$, it follows that $\eta^{n-1}(s)$ is meromorphic in $\Omega_{\frac{1}{2}}$: hence by iteration, one has $5\Rightarrow 1$ and $1\Rightarrow 5$ 
is clear.
\begin{thm} \label{th:Rie3Equ}
The following properties are equivalent:
\begin{enumerate}
\item The Riemann hypothesis holds.
\item The function $\zeta_{i}(s)$ has a meromorphic continuation in $\Omega_{\frac{1}{2}}$,  for all $i \in \kN$.
\item The function $\zeta_{i}(s)$ has an analytic continuation in $\Omega_{\frac{1}{2}}$, for all $i \in \kN$ (resp. for $i$ moving in an infinite subsequence of $\kN$).
\end{enumerate} 
\end{thm}
{\indent \it{Proof:}}
By proposition \ref{pr:EtaEtM}, proposition \ref{pr:Rie5Equ} and lemma \ref{le:ZetWexp}, one has $1\Rightarrow 2$ and $1\Rightarrow 3$. Let us assume 2). 
With notations of proposition \ref{pr:ZetZek}, one has
 $\zeta(s) \;=\; g_{i}(s)\;\zeta_{i}(s)^{i}  $ for $\kreel{1}$ and this holds by analytic continuation in 
$\Omega_{\frac{1}{2}}-
\{\text{the set of poles of}\;\zeta_{i}\}$. Let $\alpha$ be a pole of $\zeta_{i}$ in $\Omega_{\frac{1}{2}}$. Taking a small compact neighborhood $V(\alpha)$ of $\alpha$ not containing $\alpha$, we get that $\zeta_{i}(s)$ is bounded on $V(\alpha)$ hence
is analytic at $\alpha$ (ie the singularity is removable). Hence, 
$\zeta(s) \;=\; g_{i}(s)\;\zeta_{i}(s)^{i} $ for $s\in 
\Omega_{\frac{1}{2}}$ with $\zeta_{i}(s)$ analytic for all $s \in \Omega_{\frac{1}{2}}$.
Let $s_0$ be a possible root of the equation $\zeta(s)=0$ in $\Omega_{\frac{1}{2}}$. It follows that $\zeta_{i}(s)=0$ and $\zeta(s_0) \;=\; g_{i}(s)\;\zeta_{i}(s_0)^{i}$ implies that $N_{0}$ has multiplicity $i$. Since this holds for all $i \in \kN$, we get a contradiction. Hence $2\Rightarrow 1$. The proof of $3\Rightarrow 1$ is similar.
\begin{thm} \label{th:Rie3Eta}
The following properties are equivalent:
\begin{enumerate}
\item The Riemann hypothesis holds.
\item For any arithmetical list $M \subset \kP $ , the function $\frac{d}{ds}\eta_{M}(s)$ extends analytically in $\Omega_{\frac{1}{2}}$.
\item The function of $s$
$$\int_{2}^{+\infty}\; \frac{\pi(t)\ln(t)}{t^{s+1}}\;dt$$
extends as an analytic function in $\Omega_{\frac{1}{2}}$.
\end{enumerate} 
\end{thm}
{\indent \it{Proof:}}
One has $1\Rightarrow 2$ by proposition \ref{pr:Rie5Equ}. Conversely, as $\Omega_{\frac{1}{2}}$ is simply connected, one has $2\Rightarrow 1$  by lemma \ref{le:ZetWexp}.
\\
Using Abel summation, one has, for $\kreel{1}$,
$ \eta_{1}(s)\;=\; s\;\int_{2}^{+\infty}\; \pi(t)\frac{dt}{t^{s+1}}$
hence,
$$ \frac{d}{ds}\eta_{1}(s)\;=\; \int_{2}^{+\infty}\; \pi(t)\frac{dt}{t^{s+1}}
\;-\;
s\;\int_{2}^{+\infty}\; \pi(t)\frac{\ln(t)dt}{t^{s+1}}$$
Therefore, for $\kreel{1}$
$$ \;\frac{d}{ds}\eta_{1}(s)\;-\;\frac{1}{s}\eta_{1}(s)
\;=\; -s\;\int_{2}^{+\infty}\; \pi(t)\frac{\ln(t)dt}{t^{s+1}}
$$
Thus, if 3) holds, $ \;\frac{d}{ds}\eta_{1}(s)\;-\;\frac{1}{s}\eta_{1}(s)$ extends 
as an analytic function in $\Omega_{\frac{1}{2}}$, thus $\eta_{1}(s)$ has a meromorphic continuation to $\Omega_{\frac{1}{2}}$, since $\frac{d}{ds}\eta_{1}$ is always meromorphic in that set. From proposition \ref{pr:DEta}, we get $3\Rightarrow 1$. The implication $1\Rightarrow 3$ 
is similar.

\section{Summations relative to $\text{pop}_M$ }

In the previous considerations, we met subsets $\kA$ of $\kN$ having the following properties:
\begin{enumerate}
\item[1).] conditions $a\in \kA$ and $b\in \kA$ imply $ab \in \kA$
\item[2).]  conditions $a\in \kA$ and $d \;|\; a$ with  $d \in \kN$ imply $d \in \kA$
\end{enumerate}
Obviously, $1\in \kA$ when $\kA \neq \varnothing $ and:
\begin{lemma}
 There exists a unique subset $M \subset \kP$ such that
$\kA$ consists of all products (with multiplicities) of elements of $M$, ie one has $\kA=\text{pop}_M$.
\end{lemma}
We need to consider summations or formulas involving $\text{pop}_M$.
 In some cases, formulas involving $\text{pop}_M$  can be viewed as the restriction to $\text{pop}_M$ of a known formula assumed classic on $\kN$: this fact often appears for  arithmetical functions $f$  expanding some value $f(n)$ as a finite summation of others values depending on all $d \;|\;n$. Thus, proposition 7 is
the direct adaptation of Moebius formulas.

Another kind of formulas rest on the known principle of inclusion-exclusion sketched
in next lemma.
\begin{lemma}
 Let $f:\kN \rightarrow \kC$ be a function, $M \subset \kP$ with $ \kA=\text{pop}_M$. Then, 
putting $M^{*}=\kP - M$ and  $ \kB=\text{pop}_{M^{*}}$, one has:
\begin{equation} \label{principeinex}
\sum_{i \in \kA} f(i) \;=\;\sum_{k \in \kB} \mu(k) \sum_{i \in \kN} f(ki)
\end{equation}
(questions of convergence are not considered here)
\end{lemma}
Let  $M \subset \kP$ (we abbreviate $\text{pop}_M=\text{pop}$) and
$G_{pop}$ be the "geometrical" series defined for $|q|<1$ by:
\begin{equation} \label{gpop} 
G_{pop}(q) \; \equiv \; \sum_{ k\in pop}\; q^k
\end{equation}
Taking $f(i )=e^{-i t}$ for $t>0$ fixed, we get, with $ G_{pop}(e^{-t}) \equiv \sum_{n \in \text{pop} }\; e^{-nx}$ that:
$$ G_{pop}(e^{-t})  = \frac{1}{1-e^{x}} - \sum_{p \in M^{*} }\;\frac{1}{1-e^{px}}
 + \sum_{\substack{ p<q \\p,\;q \in M^{*} }}\;\frac{1}{1-e^{pqx}} -
\sum_{\substack{ p<q<r \\ p,\;q,\;r \in M^{*}}}\;\frac{1}{1-e^{pqrx}}+\cdots$$
This summation is essential in formula
\begin{equation} \label{gamzeta} 
\Gamma(s)\;Z_{M}(s) = \int_{0}^{+\infty}G_{pop}(e^{-t})\; t^s \frac{dt}{t}
\end{equation}
valid for $\kreel{1}$.
\begin{proposition} \label{pr:Lamb}
Let $M \subset \kP$  and $\text{pop}=\text{pop}_{M}$. 
\begin{enumerate}
\item The identity
$$ Z_{M}(s)\;\sum_{n \in \text{pop}}\;\frac{a_{n}}{n^s}\;=\;
\sum_{n \in \text{pop}}\;\frac{A_{n}}{n^s}$$
is equivalent to
$$ \sum_{n \in \text{pop}}\;a_{n}G_{pop}(x^n)\;=\;
\sum_{n \in \text{pop}}\;A_{n}x^n $$
\item Let 
$$ Z_{M}(s)\;\sum_{n \in \text{pop}}\;\frac{a_{n}}{n^s}\;=\;
\sum_{n \in \text{pop}}\;\frac{A_{n}}{n^s}\quad\text{and} \quad
\zeta(s)\;\sum_{n \in \text{pop}}\;\frac{a_{n}}{n^s}\;=\;
\sum_{n \in \kN}\;\frac{B_{n}}{n^s}
$$
Then, $B_{n}=A_{n}$ when $n \in \text{pop}$ and 
$B_{n}=\sum_{ i|n\; i\in \text{pop}}\; a_{i}$ when $n \notin \text{pop}$.
\end{enumerate}
\end{proposition}
Taking $\alpha_{n}=\mu(n)$ for $n\in \text{pop}$, we get, in 1):
$$ \sum_{n \in \text{pop}}\;\mu(n)G_{pop}(x^n)\;=\;x$$

Taking $\alpha_{n}=\Phi(n)$ for $n\in \text{pop}$, we get, in 1):
$$ \sum_{n \in \text{pop}}\;\Phi(n)G_{pop}(x^n)\;=\;
\sum_{n \in \text{pop}}\;n\;x^n\;=\; x \frac{d}{dx}G_{pop}(x)$$

For 2), we deduce  from Lambert series that:
$$ \sum_{n \in \text{pop}}\;a_{n}\frac{x^n}{1-x^n}\;=\;
\sum_{n \in \text{pop}}\;B_{n}x^n \;+\;
\sum_{n \in \kN-\text{pop}}\;B_{n}x^n $$

With $a_n=1$ for $n \in \text{pop}$, we get:
$$ \sum_{n \in \text{pop}}\;\frac{x^n}{1-x^n}\;=\;
\sum_{n \in \text{pop}}\; \tau(n)x^n \;+\;
\sum_{n \in \kN-\text{pop}}\;B_{n}x^n $$
where $B_{n} = \text{card}(\text{Div(n)}\cap \text{pop})$, and $\text{Div(n)}$ 
stands for the set of divisors of $n$.

In the same way,
$$ \sum_{n \in \text{pop}}\;n\;\frac{x^n}{1-x^n}\;=\;
\sum_{n \in \text{pop}}\; \sigma(n)x^n \;+\;
\sum_{n \in \kN-\text{pop}}\;B_{n}x^n $$
where $B_{n} = \sum_{ i \in \text{Div(n)}\cap \text{pop}}\;i$.

In the case $\text{pop}=\kN$, first assertion is a classical calculation met in Lambert series with $G_{pop}(x)\;=\;\frac{x}{1-x}$.
\ksaut
{\indent \it{Proof:}}
To get a visual proof, we choose $\text{pop}=\text{pop}_{2,5}=\{1,2,4,5,8,10,11,16,\cdots\}$.
Then $G_{pop}$ is the "geometrical" series
$$ G_{pop}(q) \; = \; q+q^2+q^4+q^5+q^8+q^{10}+q^{11}+\cdots$$
One has:
$$
\begin{array}{lcccccccc}
  a_1 \;G_{pop}(x)&=&a_{1}x&+\;a_{1}x^2&+\;a_{1}x^4&+\;a_{1}x^5&+\;a_{1}x^8&+\;a_{1}x^{10}& +\cdots\\
   a_2 \;G_{pop}(x^2) &= &&\;a_{2}x^2&+\;a_{2}x^4&&+\;a_{2}x^8&\;+a_{2}x^{10}&+\cdots\\
  a_4 \;G_{pop}(x^4) &= && &\;a_{4}x^4&&+\;a_{4}x^8&&+\cdots\\
  a_5 \;G_{pop}(x^5) &= &&&&\;a_{5}x^5&&+\;a_{5}x^{10}&+\cdots\\
  a_8 \;G_{pop}(x^8) &= &&&&&\;a_{8}x^8&&+\cdots\\
a_{10}\; G_{pop}(x^{10}) &= &&&&&&\;a_{10}x^{10}&+\cdots\\
\end{array}
$$
Summing columns, we get our formula since $A_{n}=0$ when $n \notin \text{pop}$ and
$A_{n}=\sum_{d|n}\;a_{d}$.\\
Assertion 2) is a simple verification.
\ksaut
An important function is, for $x>0$:
\begin{equation} \label{defnpop}
 N_{pop}(x) \; \equiv \; \text{card}([1;x]\cap \text{pop}) =\sum_{ \substack{ n\leq x \\
n\in pop}}\; 1_{pop}(n)
\end{equation}
which agree,  by Perron's formula, with
\begin{equation} \label{form:Perron}
N_{pop}(x)  = \frac{1}{2\pi i}\; \int_{\sigma-iT}^{\sigma+iT}\;Z_M(s)\frac{x^s}{s}\;ds + O(\frac{x^{\sigma+\epsilon}}{T}) 
\end{equation}
(with $\sigma>1$ and $\epsilon>0$ when $x$ is not an integer).
\begin{proposition} \label{pr:Mob}
Let $M \subset \kP$  and $\text{pop}=\text{pop}_{M}$. 
\begin{enumerate}
\item Given a function $f=f(x)$, we put
$F(x) =\;\sum_{k \in \text{pop}}\;f(kx)$. Then,
$$\sum_{n \in \text{pop}}\;a_{n}F(nx)\;=\;
\sum_{n \in \text{pop}}\;A_{n}f(nx)$$
where $a_n$ and $A_n$ are related by formula
$ \; Z_{M}(s)\;\sum_{n \in \text{pop}}\;\frac{a_{n}}{n^s}\;=\;
\sum_{n \in \text{pop}}\;\frac{A_{n}}{n^s}$.
\item Formulas
$$F(x) =\;\sum_{k \in \text{pop}}\;f(kx)$$
and
$$f(x) =\;\sum_{k \in \text{pop}}\;\mu(k)F(kx)$$
are equivalent
\end{enumerate}
\end{proposition}
{\indent \it{Proof:}}
One has
$$\sum_{n \in \text{pop}}\;a_{n}F(nx)\;=\;
\sum_{n \in \text{pop}}\;\sum_{k \in \text{pop}}\;a_{n}f(knx)=
\; \sum_{n \in \text{pop}}\;\sum_{k \;\mid \;n }\;a_{\frac{n}{k}}f(nx)$$
Since
$$ A_n \;=\; \Sigma_{ j \;\mid\; n} \; a_j $$
this proves 1), noting that
$ k \;\mid \;n $ with $n \in \text{pop}$ implies $k \in \text{pop}$.\\
For 2), taking  $a_{n}=\mu(n)$ for $n\in \text{pop}$ (zero elsewhere), one has $A_{n}=0$ for all $n$ except $A_{1}=1$, hence 
$f(x) =\;\sum_{k \in \text{pop}}\;\mu(k)F(kx)$. Conversely, from $f(x) =\;\sum_{k \in \text{pop}}\;\mu(k)F(kx)$, we obtain
$$\sum_{n \in \text{pop}}\;f(nx)\;=\;
\sum_{n \in \text{pop}}\;\sum_{k \in \text{pop}}\;\mu(k)\;F(knx)=
\; \sum_{n \in \text{pop}}\;\sum_{\substack{ k \;\mid \;n \\k \in \text{pop} }}\;\mu(\frac{n}{k})\;F(nx)$$
and our assertion follows from $\mu(1)=1$ and $\sum_{k \;\mid \;n }\;\mu(\frac{n}{k})=0$
for $n >1$ with the observation that
$$\sum_{k \;\mid \;n }\;\mu(\frac{n}{k})=\sum_{\substack{ k \;\mid \;n \\k \in \text{pop}}}
\;\mu(\frac{n}{k})$$

\begin{proposition} \label{fr:sommespop}
$\;$
\begin{enumerate}
\item For $k\in \kN$ let $G(k) =\;\sum_{i | k \;  i\in \text{pop}}\;g(i)$ where $g$ is a function defined on $pop$. Then, for $x\geq1$
\begin{equation}  \label{fr:sumonpop}
 \sum_{i \leq x \;i \in \text{pop}}\;G(i)\;=\;
\sum_{i \leq x \;i \in \text{pop}}\; N_{pop}( \frac{x}{i})\; g(i)
\end{equation}
\item  For $k\in \kN$ and $a \in \kC$, let 
$$H(k) =\;\sum_{i | k \;  i\in \text{pop}}\;h(i)\;\frac{i^a}{k^a}\qquad \text{and }
\qquad
S_{pop}(u;a)=\;\sum_{i \leq u \;  i\in \text{pop}}\;\frac{1}{i^a} $$ 
where $h$ is a function defined on $pop$. Then
\begin{equation}  \label{fr:sumHonpop}
 \sum_{i \leq x \;i \in \text{pop}}\;H(i)\;=\;
\sum_{i \leq x \;i \in \text{pop}}\; S_{pop}( \frac{x}{i};a)\; h(i)
\end{equation}
\end{enumerate}
\end{proposition}

When $ a=1$, the series  $S_{pop}(u;a)$ can be considered as  the harmonic series associated to a given pop.

{\indent \it{Proof:}}
The formula \ref{fr:sumonpop} can be considered as a particular case of \ref{fr:sumHonpop} 
with $a=0$. For \ref{fr:sumHonpop}, we expand successively $H(k)$ with $ k\in \text{pop}$ 
and $k \leq x$. The term $h(i)$ appears 
exactly one time in each expansion of $H(k)$, $k\leq x$, when $ i \;\mid\;k$: when it is it so, the coefficient of $h(i)$ is $ \frac{i^a}{k^a }$. As $\frac{k}{i} \in \text{pop}$, we see that the total coefficient of $h(i)$ is of the form 
$$\sum_{j\;\leq C(i)\; j \in \text{pop} } \; \frac{1}{j^a}$$
with $C(i)$  to be clarified .
The last term of this summation is $ \frac{i^a}{y^a }$ where $y$ is the biggest element of 
$\text{pop}\cap [1;x]$. Thus our summation consists of all $\frac{1}{j^a}$ with $j \in \text{pop}$ such that  $ ij \leq x$ with $i\in \text{pop}$. Using definition \ref{defnpop}, this
gives $C(i) = N_{pop}( \frac{x}{i})$. From
$$\sum_{j\;\leq N_{pop}( \frac{x}{i}) } \; \frac{1}{j^a} = S_{pop}( \frac{x}{i};a)$$
we get our formula \ref{fr:sumHonpop}.\\
\ksaut
\indent
In formula \ref{fr:sumonpop}, we choose for $g$ the Euler $\Phi$-function. Since 
$\sum_{i | k }\;\Phi(i)=k$, we get
\begin{equation}  \label{fr:phipop}
 \sum_{i \leq x \;i \in \text{pop}}\;i\;=\;
\sum_{i \leq x \;i \in \text{pop}}\; N_{pop}( \frac{x}{i})\; \Phi(i)
\end{equation}

Taking for $g$ the function $g(i)=\mid\mu(i)\mid$, we get
\begin{equation}  \label{fr:2powipop}
 \sum_{i \leq x \;i \in \text{pop}}\;2^{\nu(i)}\;=\;
\sum_{i \leq x \;i \in \text{pop}}\; N_{pop}( \frac{x}{i})\; \mid\mu(i)\mid
\end{equation}
where $\nu(n)$ is the number of distinct prime factors of $k$: this follows from
$\sum_{i | k} \;\mid\mu(i)\mid = 2^{\nu(k)}$. Thus
$$ 
 \sum_{i \leq x \;i \in \text{pop}}\;2^{\nu(i)}\;=\;N_{pop}(x)+
\sum_{p \in \text{pop}}\; N_{pop}( \frac{x}{p}) 
+ \sum_{\substack{p<q \\ \;p \in \text{pop}  \;q \in \text{pop}}}\; N_{pop}( \frac{x}{pq})
+ \sum_{\substack{p<q<r \\ \;p \in \text{pop}  \;q \in \text{pop}\\\;r \in \text{pop}}}\; N_{pop}( \frac{x}{pqr})\;+\cdots
$$

If we choose for $g$ the Moebius $\mu$-function, from $G(i)=\sum_{i | k }\;\mu(i)=0$ for $i>1$ and $G(1)=1$, we get for $x>1$
\begin{equation}  \label{fr:muGpop}
 \sum_{i \leq x \;i \in \text{pop}}\; N_{pop}( \frac{x}{i})\; \mu(i)=1
\end{equation}

Taking for $g$ the function $g(i)=\frac{\mu(i)}{i}$, from formula $\sum_{i | k} \;\frac{\mu(i)}{i} = 
\frac{\Phi(k) }{k}$, we get
\begin{equation}  \label{fr:phisuripop}
 \sum_{i \leq x \;i \in \text{pop}}\;\frac{\Phi(i)}{i}\;=\;
\sum_{i \leq x \;i \in \text{pop}}\; N_{pop}( \frac{x}{i})\;\frac{\mu(i)}{i} 
\end{equation}

Taking for $g$ the function $g(i)=1$ when $i \in \text{pop}$, from formula $\sum_{i | k}
\;1=\tau(k)$ (the number of positive divisors of $k$), we get
\begin{equation}  \label{fr:phidivpop}
 \sum_{i \leq x \;i \in \text{pop}}\;\tau(i)\;=\;
\sum_{i \leq x \;i \in \text{pop}}\; N_{pop}( \frac{x}{i}) 
\end{equation}

Replacing $g$ by the function $g(i)=1$ when $i \in \text{pop}$ and $i$ is a prime number, we get
\begin{equation}  \label{fr:phinupop}
 \sum_{i \leq x \;i \in \text{pop}}\;\nu(i)\;=\;
\sum_{\substack{i \leq x \;i \in \text{pop}\\ i\;  \text{prime}}}\; N_{pop}
( \frac{x}{i}) 
\end{equation}

Taking for $g$ the function $g(i)=i^a$ when $i \in \text{pop}$, from formula $\sum_{i | k}
\;i^a=\sigma_{a}(k)$ (the sum of powers of positive divisors of $k$), we get
\begin{equation}  \label{fr:phisigapop}
 \sum_{i \leq x \;i \in \text{pop}}\;\sigma_{a}(i)\;=\;
\sum_{i \leq x \;i \in \text{pop}}\; N_{pop}( \frac{x}{i})\;i^a 
\end{equation}

Noting that $\sum_{i | k}
\;i^{-a}=\frac{1}{k^a}\sigma_{a}(k)$, we get for $g(i)=i^{-a}$ 
\begin{equation}  \label{fr:phisigminapop}
 \sum_{i \leq x \;i \in \text{pop}}\;\frac{\sigma_{a}(i)}{i^a}\;=\;
\sum_{i \leq x \;i \in \text{pop}}\; N_{pop}( \frac{x}{i})\;\frac{1}{i^a} 
\end{equation}

We now use formula \ref{fr:sumHonpop}. With $h(i)= 1$ for all $i$, one has $H(k)= \frac{\sigma_a(k)}{k^a}$, hence
\begin{equation}  \label{fr:sigaHonpop}
 \sum_{i \leq x \;i \in \text{pop}}\;\frac{\sigma_a(i)}{i^a}\;=\;
\sum_{i \leq x \;i \in \text{pop}}\; S_{pop}( \frac{x}{i};a)
\end{equation}

With $h(i)= \frac{1}{i^a}$ for all $i$, one has $H(k)= \frac{\tau(k)}{k^a}$, hence
\begin{equation}  \label{fr:tauaHonpop}
 \sum_{i \leq x \;i \in \text{pop}}\;\frac{\tau(i)}{i^a}\;=\;
\sum_{i \leq x \;i \in \text{pop}}\; S_{pop}( \frac{x}{i};a)\;\frac{1}{i^a}
\end{equation}

With $h(i)= \frac{1}{i^a}$ when $i$ is a prime number and $0$ elsewere, one has $H(k)= \frac{\nu(k)}{k^a}$, hence
\begin{equation}  \label{fr:nuaHonpop}
 \sum_{i \leq x \;i \in \text{pop}}\;\frac{\nu(i)}{i^a}\;=\;
\sum_{\substack{i \leq x \;i \in \text{pop}\\ i\;  \text{prime}}}\; S_{pop}( \frac{x}{i};a)\;\frac{1}{i^a}
\end{equation}

More generally, for $k\in \kN$, let 
$F(k) =\;\sum_{i | k \;  i\in \text{pop}}\;f(i)g(\frac{k}{i})$
where $h$ and $g$ are functions defined on pop. Then, the  same calculation yields
\begin{equation}  \label{fr:sumHKpop}
 \sum_{i \leq x \;i \in \text{pop}}\;F(i)\;=\;
\sum_{i \leq x \;i \in \text{pop}}\; f(i) 
\sum_{j \leq \frac{x}{i} \;j \in \text{pop}}\;g(j)
\end{equation}

Thus, with $x$ infinite, $f(i)=a^i$, $g(i)=b^i$, we obtain, for $\mid a \mid <1$ and 
 $\mid b \mid <1$
\begin{equation} \label{fr:gpopgpop}
G_{pop}(a)\;G_{pop}(b) = \sum_{n \in \text{pop}\;n \geq 1}\big(
\sum_{d | n }a^d b^{\frac{n}{d}} \big)
\end{equation}

One also has the formula
\begin{equation}  \label{fr:sumHKdesipop}
 \sum_{i \leq x \;i \in \text{pop}}\;i\;=\;
\sum_{i \leq x \;i \in \text{pop}}\; \mu(i) 
\sum_{j \leq \frac{x}{i} \;j \in \text{pop}}\;\sigma(j)
\;=\;
\sum_{i \leq x \;i \in \text{pop}}\; \sigma(i) 
\sum_{j \leq \frac{x}{i} \;j \in \text{pop}}\;\mu(j)
\end{equation}
due to $\sum_{i | k }\;\sigma(i)\mu(\frac{k}{i} ) = k$ for $k \in \kN$. And
$\sum_{i | k }\;\tau(i)\mu(\frac{k}{i} ) = 1$ for $i \in \kN$ yields
\begin{equation}  \label{fr:sumHKNpop}
 N_{pop}(x) \;=\;\sum_{i \leq x \;i \in \text{pop}}\; \mu(i) 
\sum_{j \leq \frac{x}{i} \;j \in \text{pop}}\;\tau(j)
\;=\;
\sum_{i \leq x \;i \in \text{pop}}\; \tau(i) 
\sum_{j \leq \frac{x}{i} \;j \in \text{pop}}\;\mu(j)
\end{equation}
$\tau(n)$ being the divisor function.

Let $\kA$ be a subset of $\kN$. We put ${\kA }^{+} =  \bigcup_{i \in \kA} [i;i+1[ = \kA+[0;1[$
and, for $ n\in \kN$,
${\kA }_{n}^{+} =  \bigcup_{i \in \kA \; i\leq n} [i;i+1[ $.\\
The comparison of a series with an integral takes here the following aspect
\begin{proposition}
Let $\kA$ be a subset of $\kN$ and $f:[1;+\infty[\longrightarrow \kR$ be a continuous positive decreasing function on the interval $[1;+\infty[$ such that $\text{limit}_{t \rightarrow+\infty}f(t)=0$.\\
Then, the sequence:
$$ u_{n} =  \sum_{
		i \leq n 
		\; i \in \kA} f(i) - \int_{{\kA}_{n}^{+}} f(t)dt $$
has a finite limit $C_{A}$ when $n \rightarrow +\infty$.
\end{proposition}
{\it{Proof:}}
We modify the function $f$ on interval $[i;i+1[$ whenever $i \notin \kA$ 
making it affine on such intervals (values $f(i)$ and $f(i+1)$ being unalterated).
Let $h:[1;+\infty[\longrightarrow \kR$ be the continuous positive decreasing function so obtained.
It is known that the sequence:
$$
\begin{aligned}   v_{n} &=  \sum_{i=1}^{i=n} h(i) - \int_{1}^{n+1} h(t)dt\\
 &= h(1) - \int_{1}^{2} h(t)dt  +\cdots
+h(i) - \int_{i}^{i+1} h(t)dt +\cdots + h(n) - \int_{n}^{n+1} h(t)dt
 \end{aligned}
$$
has a finite limit  when $n \rightarrow +\infty$. Whenever $i \notin \kA$, one has:
$$ h(i) - \int_{i}^{i+1} h(t)dt \;=\; \frac{1}{2} \big(f(i)-f(i+1)\big)$$
thus
$$ u_{n} = v_{n} - \frac{1}{2}\;\sum_{i \notin \kA \;i \leq n}  \big(f(i)-f(i+1)\big)$$
Our conclusion follows since $\sum_{i \notin \kA}  \big(f(i)-f(i+1)\big)$ satisfies  hypothesis of the theorem of alternating convergent series.

\begin{proposition}Let $\kA$ be a subset of $\kN$ and $Z(s)$ be the function 
$ Z(s)=\Sigma_{n \in \kA}\;\frac{1}{n^{s}}$ for $\kreel{1}$. One has, for $\kreel{1}$:
$$ Z(s) = \int_{{\kA }^{+}} \frac{dt}{t^{s}} \;+\; h(s)$$
where $h(s)$ is analytic in the half-plane $\kreel{0}$.
\end{proposition}
{\it{Proof:}}
Classic demonstration (see [Nar] p.209 when $\kA=\kN$) adapts itself without any trouble. One has, for $n \in \kN$ and $\kreel{1}$:
$$\int_{0}^{1} \frac{dt}{(t+n)^{s}}\;=\;\frac{1}{(1+n)^{s}}\;+\; s\;\int_{0}^{1} \frac{t\;dt}{(t+n)^{s+1}}$$
Summing for $n \in \kA$, we get
$$ Z(s)\;=\;\sum_{n \in \kA}\bigg(\frac{1}{n^{s}}-\frac{1}{(1+n)^{s}}\bigg)
- \sum_{n \in \kA}\;s\;\int_{0}^{1} \frac{t\;dt}{(t+n)^{s+1}}
+ \sum_{n \in \kA} \int_{0}^{1} \frac{dt}{(t+n)^{s}}$$
The first expression of the right-hand side is an alternate Dirichlet series clearly 
analytic for $\kreel{0}$.
We define, for $t \in [0;1]$ : 
$$g_{\kA}(t)=\sum_{n \in \kA}\frac{t}{(t+n)^{s+1}}$$
The series of $g_{\kA}(t)$ converges uniformly for $t \in [0;1]$ and $s$  in any given compact of $\kreel{0}$, hence 
$$\sum_{n \in \kA}\;s\;\int_{0}^{1} \frac{t\;dt}{(t+n)^{s+1}}\;=\;s \int_{0}^{1}g_{\kA}(t)dt$$
defines an analytic function in that half-plane.\\ 
Since
$$ \sum_{n \in \kA} \int_{0}^{1} \frac{dt}{(t+n)^{s}} = \int_{{\kA }^{+}} \frac{dt}{t^{s}}$$
this proves the lemma.
\ksaut
\indent 
One has, for $\kreel{1}$:
\begin{equation} \label{zetanpop} 
Z_{M}(s) = \sum_{n \in \kN} N_{pop_{M}}(n)\bigg[ \frac{1}{n^s}-\frac{1}{(n+1)^s}\bigg]
=s\;\int_{1}^{+\infty}N_{pop_{M}}(t)\frac{dt}{t^{s+1}}
\end{equation}
and, by Abel summation, for $t>0$:
$$
G_{pop}(e^{-t}) \;=\; \sum_{k=1}^{+\infty}N_{pop}(k)\; t \int_{k}^{k+1}e^{-tu}\;du \;=\;
t \int_{1}^{+\infty}N_{pop}(u) \;e^{-tu}\;du $$
Now, from formula \ref{principeinex} with $g(i)=1$ for all $i\leq n$ and $0$ elsewere, in the case $\text{pop}_M=\text{pop}_{2,5}$, we obtain
$$ N_{pop}(n)= n - \bigg[ \frac{n}{3}\bigg]- \bigg[ \frac{n}{7}\bigg]-\cdots
+ \bigg[ \frac{n}{21}\bigg]+\cdots$$
which is close to
$$ n(1-\frac{1}{3})(1-\frac{1}{7})(1-\frac{1}{13})\cdots \approx \;n\;\prod_{p \leq n \; p \in \text{pop}_{3,7}}(1-\frac{1}{p})\;
\sim \; \alpha_{3,7}\; \frac{n}{\sqrt{\ln(n)} }$$
for some constant $\alpha_{3,7}$ (by proposition below). Numerical computations directly checked on $\text{pop}_{2,5} = \{1,2,4,5,8,10,11,16,17,\cdots\}$ gives $\alpha_{3,7}\sim 0,736\cdots$ with
an oscillating limit.\\
More generally, it is highly probable that, for an arithmetical list $M$ having a reason $r>0$, one has, for some suitable constant $A_M$:
$$ N_{pop}(x)\; \sim \; A_M  \;\;x\; \frac{\sqrt[r]{\ln(x)}}{\ln(x)}\quad \text{as}\quad  x\rightarrow +\infty$$
This property is related to formula $Z_{M}(s) \sim \frac{\text{Cte}}{(s-1)^{\frac{1}{r}}}$ associated with Perron's formula \ref{form:Perron} or a possible improvement of Ikehara-Wiener theorem.\\
Mertens's formulas takes the form:
\begin{proposition} \label{pr:Mert}
Let $M$ be an arithmetical list of $\kP$,  $r>0$ being its reason. There exists a constant $\gamma_M$ such that one has:
\begin{equation} \label{fr:mertens1}
\prod_{\substack{p \leq x \\ p \in M}}(1-\frac{1}{p})
= \;\frac{e^{-\gamma_M}}{\sqrt[r]{\ln(x)}} + O\bigg( \frac{1}{\ln(x)\sqrt[r]{\ln(x)}}\bigg)
\end{equation}
\begin{equation} \label{fr:mertens2}
\frac{1}{
\prod_{\substack{p \leq x \\ p \in M}}(1-\frac{1}{p})
}
= \; e^{\gamma_M} \sqrt[r]{\ln(x)} + O\bigg( \frac{\sqrt[r]{\ln(x)} }{\ln(x) }\bigg)
\end{equation}
\end{proposition}

When $M=\kP$, one has $r=1$ and $\gamma_M$ is the Euler constant $\gamma$.
\ksaut 
{\indent \it{Proof:}}
Demonstration proposed here is inspired, in the main lines, by that of [Nat] p. 278. One has  
$$ \sum_{p \leq x\;,\;p\in M} \frac{1}{p} = \sum_{p \leq x\;,\;p\in M} \frac{\ln(p)}{p}\;\frac{1}{\ln(p)}
=\sum_{2 \leq n \leq x} f(n)g(n)
$$
with $g(t) = \frac{1}{\ln(t)}$ for $t>2$, $f(n)=\frac{\ln(p)}{p} $ if $ n=p \in M$ and $f(n)=0$ otherwise.\\
From known formula
\begin{equation}
 \sum_{p \leq x}\;\frac{\ln(p)}{p} = \ln(x) + O(1) 
\end{equation}
we first deduce that
\begin{equation} \label{fr:lnpsurp}
\sum_{p \leq x\;,\;p\in M} \frac{\ln(p)}{p} = \frac{1}{r}\ln(x) + O(1) 
\end{equation}
Indeed, if 
$$ p_1 \leq q_1 \leq p_2 \leq q_2 \leq \cdots\leq p_n \leq q_n \leq p_{n+1} \leq \cdots $$
are real numbers tending to infinity as $n \rightarrow +\infty$, one must have
$$ \sum_{p_i \leq x}\;f(p_i) - \sum_{q_i \leq x}\;f(q_i)\;=\; O(1) $$
since
$$ f(p_1) - f(q_1) + f(p_2) - f(q_2) +\cdots +f(p_n) - f(q_n) + \cdots $$
is an alternating convergent sequence (and partial summations gives rise to two monotonic sequences having the same limit): thus, we adapt the proof of proposition \ref{pr:EtaEtM}  taking $s=1$, $w_j(s)=O(1)$ and formula \ref{fr:lnpsurp} follows.\\
Let
$$ F(t) = \sum_{p \leq t\;,\;p\in M} \frac{\ln(p)}{p} =  \sum_{p \leq t\;,\;p\in M} f(n)$$
Then $F(t) = \frac{1}{r}\ln(t) + r(t)$ with $r(t)=O(1)$. Our next goal is to get formula 
\begin{equation} \label{fr:onesurp}
\sum_{p \leq x\;,\;p\in M} \frac{1}{p} = \frac{1}{r}\ln(\ln(x)) + b + O\bigg(\frac{1}{\ln(x)}
\bigg) 
\end{equation}
for some fixed $b \in \kR$. This is a consequence of partial summation since
$$  \int_{x}^{+\infty}\; \frac{r(t)}{t (\ln(t))^2}\; dt = O\bigg(\frac{1}{\ln(x)}\bigg)$$
and
$$
\begin{aligned}
\sum_{2 \leq n \leq x} f(n)g(n)
&= F(x)g(x) - \int_{2}^{+\infty}\;F(t)g'(t)dt \\
 &= \frac{1}{r} + O\bigg(\frac{1}{\ln(x)}\bigg)
+ \frac{1}{r} \int_{2}^{x}\; \frac{1}{t \ln(t)}\; dt
+ \frac{1}{r} \int_{2}^{x}\; \frac{r(t)}{t (\ln(t))^2}\; dt\\
&=\frac{1}{r} + \frac{1}{r}\ln(\ln(x))-\frac{1}{r}\ln(\ln(2)) + \frac{1}{r}\int_{2}^{+\infty}\; \frac{r(t)}{t (\ln(t))^2}\; dt + O\bigg(\frac{1}{\ln(x)}\bigg)
\end{aligned}
$$
leading to formula \ref{fr:onesurp} with 
$$ b = \frac{1}{r} -\frac{1}{r}\ln(\ln(2)) + \frac{1}{r}\int_{2}^{+\infty}\; \frac{r(t)}{t (\ln(t))^2}\; dt$$
It follows that
$$ \ln(\prod_{\substack{p \leq x \\ p \in M}}(1-\frac{1}{p})^{-1}
)=  -\sum_{p \leq x\;,\;p\in M} \ln(1 - \frac{1}{p})
= \sum_{p \leq x\;,\;p\in M} \frac{1}{p}  + \sum_{p\in M} \sum_{k=2}^{\infty}\frac{1}{kp^k}-
\sum_{p\in M\;p>x} \sum_{k=2}^{\infty}\frac{1}{kp^k}$$
Noting that
$$ \sum_{p\in M} \sum_{k=2}^{\infty}\frac{1}{kp^k} \leq 
\sum_{p\in \kP} \sum_{k=2}^{\infty}\frac{1}{kp^k} < +\infty
$$
$$\sum_{p\in M\;p>x} \sum_{k=2}^{\infty}\frac{1}{kp^k} \leq
\sum_{p\in \kP\;p>x} \sum_{k=2}^{\infty}\frac{1}{kp^k} \leq \frac{2}{x}
$$
we find a  constant $\gamma_M$ such that:
$$ \ln(\prod_{\substack{p \leq x \\ p \in M}}(1-\frac{1}{p})^{-1}= \frac{1}{r}\text{ln(ln(}x)) + \gamma_{M} + O\bigg( \frac{1}{\ln(x)}\bigg)$$
Thus,
$$
\begin{aligned}
\frac{1}{
\prod_{\substack{p \leq x \\ p \in M}}(1-\frac{1}{p})
}
&= \; e^{\gamma_M} \sqrt[r]{\ln(x)}\; e^{O( \frac{1}{\ln(x)}) }\\
&= \; e^{\gamma_M} \sqrt[r]{\ln(x)}\;\bigg(1 + O\bigg( \frac{1}{\ln(x)}\bigg)\bigg)
\end{aligned}
$$
which is formula \ref{fr:mertens2}. Taking the inverse, we get formula \ref{fr:mertens1}.\\ 
\begin{proposition} \label{pr:Mesch}
Let $M$ be an arithmetical list of $\kP$ with a reason $r>0$. We assume that, for large $t$:
\begin{equation} \label{fr:equilog}
 N_{pop}(t)\; \sim \; A_M  \;t \frac{\sqrt[r]{\ln(t)}}{\ln(t)}
\end{equation}
for some suitable constant $A_M$. Then, one has:
$$ \sum_{n \in \text{pop}_M \;,\; n \leq x}\; \frac{1}{n}\quad \sim \quad r\;A_{M} \;\sqrt[r]{\ln(x)} \quad \text{as} \quad x \rightarrow +\infty$$
\end{proposition}
{\indent \it{Proof:}}
With Abel summation, putting $ m = [x]$, one has:
$$
 \sum_{n \in \text{pop}_M \;,\; n \leq x}\;\frac{1}{n}\;
=\;\int_{1}^{m}\;N_{pop}(t)\;\frac{dt}{t^2}+ \frac{1}{m} N_{pop}(m)
$$ 
Clearly, $ \frac{1}{m} N_{pop}(m)
\; \sim \; A_M  \; \frac{\sqrt[r]{\ln(m)}}{\ln(m)}\quad\rightarrow\quad 0
 \quad \text{as} \quad m \rightarrow +\infty$.
For large $A$, one may write, for $ t \geq A$, 
$ N_{pop}(t)\;=\; (A_M+o(1))  \;t\; \frac{\sqrt[r]{\ln(t)}}{\ln(t)}$, hence for $ m \geq A$:
$$\int_{A}^{m}\;N_{pop}(t)\;\frac{dt}{t^2}=
 (A_M+o(1))  \; \int_{A}^{m}\; \frac{\sqrt[r]{\ln(t)}}{\ln(t)} \;\frac{dt}{t}= (A_M+o(1))\;\bigg[r\;\sqrt[r]{\ln(t)}\bigg]_{t=A}^{t=m}$$
$$
= r\;(A_M+o(1))\;\big( \sqrt[r]{\ln(m)}-  \sqrt[r]{\ln(A)}\big)\; \sim \;
r\;A_M\;\sqrt[r]{\ln(m)}\quad \text{as} \quad m \rightarrow +\infty
$$
Now, with $ A$ fixed, as $ 1 \leq N_{pop}(t)\leq t$ for any $ t \geq1$, one has:
$$ 0 \; \leq \; \frac{\int_{1}^{A}\;N_{pop}(t)\;\frac{dt}{t^2}}{r\;A_M\;\sqrt[r]{\ln(m)}}
\;\leq\; \frac{\ln(A)}{r\;A_M\;\sqrt[r]{\ln(m)}}\quad  \rightarrow \quad 0\quad
\text{when} \quad m \rightarrow +\infty
$$
Hence, $$ \int_{1}^{m}\;N_{pop}(t)\;\frac{dt}{t^2}\; \sim \;
r\;A_M\;\sqrt[r]{\ln(m)}\quad \text{as} \quad m \rightarrow +\infty
$$
This proves the proposition.
\begin{proposition} \label{pr:SigExp}
Let $M$ be an arithmetical list of $\kP$ with a reason $r>0$. We assume that 
$$ N_{pop}(t)\; \sim \; A_M  \;\;t\; \frac{\sqrt[r]{\ln(t)}}{\ln(t)}\quad \text{as}\quad  t\rightarrow +\infty$$
 Then, one has,  for $x > 0$:
$$ \sum_{n \in \text{pop}_M }\; e^{-nx}  
\;=\; 
x\;\int_{1}^{+\infty}\; N_{pop}(t)\;e^{-tx}dt
$$
For $x$ near 0, one has
$$
\sum_{n \in \text{pop}_M }\; e^{-nx}  
\;\sim\;
\frac{A_M }{x } \; 
\frac{\sqrt[r]{\ln(\frac{1}{x} )}}{\ln(\frac{1}{x})}
$$
\end{proposition}
In formula \ref{gamzeta}, using $\int_{0}^{+\infty}=\int_{0}^{1}+\int_{1}^{+\infty}$, we obtain  that the function 
$$ \varphi:\; s \longrightarrow \int_{1}^{+\infty}\sum_{n \in \text{pop}_M }\; e^{-nt} \; t^s \frac{dt}{t} $$
 clearly defines an analytic function $\varphi(s)$ for $\kreel{0}$.  The remaining integral can be written:
$$\int_{0}^{1}\sum_{n \in \text{pop}_M }\; e^{-nt} \; t^s \frac{dt}{t} 
= \int_{0}^{1} \frac{A_M }{t } \; 
\frac{\sqrt[r]{\ln(\frac{1}{t} )}}{\ln(\frac{1}{t})} \; t^s \frac{dt}{t} +
\int_{0}^{1} h(t)\;t^s \frac{dt}{t}
$$
for some unknown function $h(t)$, weak near $t=0$. Since one has, for $\kreel{1}$
$$ \int_{0}^{1}\;\frac{1 }{t }\;\frac{\sqrt[r]{\ln(\frac{1}{t} )}}{\ln(\frac{1}{t})}\;t^{s-1}dt \;=\; 
\frac{1}{\sqrt[r]{s-1}}\;\Gamma(\frac{1}{r} )$$
we get that:
$$ \Gamma(s)\;Z_{M}(s) =   
\frac{A_M}{\sqrt[r]{s-1}}\;\Gamma(\frac{1}{r}) + \int_{0}^{1} h(t)\;t^s \frac{dt}{t}
+\varphi(s)$$
which  is consistent with the property \ref{form:akzetak}. 
\ksaut
\indent It is also known that $\sum_{n \in \kN \; n\geq 1 }\; e^{-nx} 
\;=\; \frac{1}{x} -\frac{1}{2} + \sum_{n=1}^{+\infty}\; \frac{2x}{x^2+4 n^2 \pi^2} $ and important features of $\zeta(s)$ are obtained from this formula. When $M$ is an arithmetical list with a reason $r$, 
simple arguments seems to indicate that the corresponding development of
$ \sum_{n \in \text{pop}_M }\; e^{-nx}$  
looks like
$$\sum_{n \in \text{pop}_M }\; e^{-nx} 
\;=\; \frac{B_r }{x } \; 
\frac{\sqrt[r]{\ln(\frac{1}{x} )}}{\ln(\frac{1}{x})} + \sum_{n \in \text{pop}_M }\; \frac{2x}{x^2+4 n^2 \pi^2} + A(x)$$
where $A(x)$ is a new expression to be found. A complete knowledge of  this expansion may be essential to the solution of Riemann's conjecture.
\ksaut
{\indent \it{Proof of proposition \ref{pr:SigExp}
:}}
Let $m \in \kN,\; m \geq 2$. By Abel summation, 
$$
\sum_{1\leq n \leq m\; n \in \text{pop}_M }\; e^{-nx} \;
=\;x\; \int_{1}^{m}\;N_{pop}(t)\;e^{-tx}dt+ e^{-mx} N_{pop}(m)
$$
Clearly, for $ x>0$, 
$$
e^{-mx} N_{pop}(m)
\; \sim \; m\; \frac{A_M\;\sqrt[r]{\ln(m)}}{\ln(m)}e^{-mx}\;
\longrightarrow 0 \qquad \text{when} \qquad  m \rightarrow +\infty$$ 
As $ 1 \leq N_{pop}(t)\leq t$ for any $ t \geq1$, the function $t \rightarrow N_{pop}(t)\;e^{-tx} $ is integrable on the interval $ [1;+\infty[$, hence letting $ m \rightarrow +\infty$,
one has, for fixed $x>0$
$$
\sum_{1\leq n \; n \in \text{pop} }\; e^{-nx} \;
=\;x\; \int_{1}^{+\infty}\;N_{pop}(t)\;e^{-tx}dt
$$
since $\sum_{1\leq n \; n \in \text{pop} }\; e^{-nx}$ is a convergent series.\\
Let $A>1$. We put $f(t) = t\;\frac{\sqrt[r]{\ln(t)}}{\ln(t)} $ for $t\geq A$. 
One has, for large $t$: 
$$ \frac{f'(t) }{ f(t)} = \frac{1}{t} + (\frac{1}{r}-1)\;\frac{1}{t\ln(t)}
\;\sim\;\frac{1}{t} $$
and it follows from Laplace's method that, for $x$ near 0:
$$
\int_{A}^{+\infty}\; f(t)\;e^{-tx}dt
 \;\sim\; \frac{\Gamma(2) }{x } \; f(\frac{1 }{x })= \frac{1 }{x^2 } \; 
\frac{\sqrt[r]{\ln(\frac{1}{x} )}}{\ln(\frac{1}{x})}
$$
For any $ 0<\epsilon <1$, one can find $A=A(\epsilon)$ such that, for $t\geq A$, one has
$$ (1-\epsilon)\;A_M  \; f(t)
 \;\leq \; N_{pop}(t) \;\leq \;(1+\epsilon)\;A_M  \; f(t)
 $$
leading to
 $$ 
\begin{array}{ll}
 (1-\epsilon)\;A_M \;x\;\int_{A}^{+\infty}\;f(t)\;e^{-tx}dt
 \;&\leq \;x\; \int_{A}^{+\infty}\;N_{pop}(t)\;e^{-tx}dt 
\\
\;&\leq \;
(1+\epsilon)\;A_M \;x\;\int_{A}^{+\infty}\;f(t)\;e^{-tx}dt
\end{array}
$$
for any $x>0$, thus it remains to show that
$$  \frac{x \; \int_{1}^{A}\;N_{pop}(t)\;e^{-tx}dt }
{\frac{A_M }{x } \; 
\frac{\sqrt[r]{\ln(\frac{1}{x} )}}{\ln(\frac{1}{x})}
 } \quad  \rightarrow 0 \quad
\text{as} \quad x \rightarrow 0
$$
for $A$ fixed, 
which easily follows from $\int_{1}^{A}\;N_{pop}(t)\;e^{-tx}dt \leq A^2$.
This proves the proposition.\\
\begin{center}
{\bf{References}}
\end{center}
$[Gr$-$Sc]$ \ E.Grosswald, F.J.Schnitzler. A class of modified $\zeta$ and $L$-functions. Pacific J.Math., 74, 357-364, 1978.\\
$[Iv]$ \ A.Ivic. Mean values of the Riemann Zeta-function. Tata Institute of Fundamental Research. Springer, Berlin-Heidelberg 1991.\\
$[Nar]$ \ W.Narkiewicz. The development of Prime Number Theory. Springer, Berlin-Heidelberg 2000.\\
$[Nat]$ \ M.B.Nathanson. Elementary Methods in Number Theory. volume 195 of Graduate Texts in Mathematics. Springer-Verlag, New York Berlin Heidelberg 2000.\\
$[Pat]$ \ S.Patterson. An introduction to the theory of the Riemann zeta function. Cambridge Studies in Advanced Mathemaics {\bf{14}}. Cambridge University Press, 1988.\\
$[Ser]$ \ J-P.Serre. Cours d'arithm\'etique. P.U.F Paris, 1970. \\
$[Tit]$ \ E.C.Titchmarsh. The Theory of Riemann Zeta-function. Oxford.  Clarendon Press, 1951.\\

\end{document}